\documentclass[a4paper]{article}
\usepackage[a4paper, total={6in, 8in}]{geometry}
\usepackage[english]{babel}
\usepackage[utf8]{inputenc}
\usepackage{amsmath}
\usepackage{amssymb}
\usepackage{graphicx}
\usepackage{hyperref}
\usepackage[colorinlistoftodos]{todonotes}

\usepackage[
    backend=biber,
    style=numeric-comp,
    sorting=none,
  ]{biblatex} 
\usepackage{authblk}  

\title{A Novel Route to Oscillations via non-central SNICeroclinic Bifurcation: Unfolding the Separatrix Loop Between a Saddle-Node and a Saddle}

\author[1,2]{Kateryna Nechyporenko}
\author[2]{Peter Ashwin}
\author[1,2]{Krasimira Tsaneva-Atanasova}

\affil[1]{Living Systems Institute, University of Exeter, Stocker Road, Exeter EX4 4PY, UK}
\affil[2]{Department of Mathematics and Statistics, University of Exeter, Harrison Building, Exeter EX4 4QF, UK}

\date{}

\begin{document}

\maketitle

\begin{abstract}
    In this paper, we investigate saddle–node to saddle separatrix-loops that we term SNICeroclinic bifurcations. They are generic codimension-two bifurcations involving a heteroclinic loop between one non-hyperbolic and one hyperbolic saddle. A particular codimension-three case is the non-central SNICeroclinic bifurcation. We unfold this bifurcation in the minimal dimension (planar) case where the non-hyperbolic point is assumed to undergo a saddle-node bifurcation. Applying the method of Poincar\'{e} return maps, we present a minimal set of perturbations that captures all qualitatively distinct behaviors near a non-central SNICeroclinic loop. Specifically, we study how variation of the three unfolding parameters leads to transitions from heteroclinic and homoclinic loops; saddle-node on an invariant circle (SNIC); and periodic orbits as well as equilibria. We show that although the bifurcation has been largely unexplored in applications, it can act as an organizing center for transitions between various types of saddle-node and saddle separatrix loops. It is also a generic route to oscillations that are both born and destroyed via global bifurcations, compared to the commonly observed scenarios involving local (Hopf) and in some cases global (homoclinic or SNIC) bifurcations.  
\end{abstract}

\textbf{Keywords:} Global bifurcations driving oscillatory patterns, Saddle separatrix-loops, Mechanisms of self-sustained oscillations; Dynamical systems in interdisciplinary research, Non-hyperbolic heteroclinic loop

\textbf{AMS:} 92B25, 34C15, 34C23, 34C37


\section{Introduction}

Saddle-node separatrix loops involving homoclinic and heteroclinic connections are important for understanding dynamical systems and their applications. For example, homoclinic bifurcations have been implicated in the control of different types of dynamic behavior including periodic solutions in climate models \cite{shilnikov_bifurcation_1995, osullivan_rate-induced_2023}, lasers \cite{shilnikov_homoclinic_1997,wieczorek_bifurcations_2005, krauskopf_excitability_2003} and biochemical models \cite{kaas-petersen_homoclinic_1988,kazmierczak_homoclinic_2002} as well as neurodynamics \cite{izhikevich_neural_2000, barrio_dynamics_2024}.
Meanwhile, heteroclinic dynamics has been implicated in studies of information flow in cognitive processes \cite{aravind_relaxation_2023,rabinovich_information_2012}, decision making \cite{varona_hierarchical_2016-1}, memory \cite{afraimovich_sequential_2015} and attention \cite{rabinovich_hierarchical_2015}. However, most of the applications (especially of heteroclinic dynamics) have examined the bifurcation problems of orbits that are connected via hyperbolic equilibria, while the case of non-hyperbolic equilibria has received significantly less attention.

In particular, there has been a body of previous work focusing on the saddle-node homoclinic loop (SNHL) bifurcation (or saddle-node separatrix loop), a codimension-2 bifurcation of two-dimensional vector fields which represents a saddle-node whose separatrix forms a closed curve in the boundary of the two-dimensional invariant manifold \cite{guckenheimer_multiple_1986}. This bifurcation has been introduced for the first time in the work of Luk'yanov \cite{lukyanov_bifurcations_1982} followed by Schecter \cite{schecter_saddle-node_1987}, who performed an analytical study of the planar case using Mel'nikov functions.
The results were generalized for any finitely-dimensional system by Chow and Lin in \cite{chow_bifurcation_1990} applying a variety of techniques (exponential dichotomy; Mel'nikov function; smooth foliation; and Shil’nikov’s central ideal). Their work appeared in parallel with further analysis by Deng \cite{deng_homoclinic_1989} aiming to unify the methods for studying homoclinic bifurcations.

Non-central saddle-node homoclinic (NCH) (another name found in literature describing SNHL) bifurcations can be found in biochemical systems, such as a model of GTPase activation \cite{zmurchok_local_2023}, chemical reaction-diffusion system \cite{dilao_excitability_2004}, and mathematical model of biochemical reactions \cite{borisuk_bifurcation_2008}. Maruyama \textit{et al.} \cite{maruyama_analysis_2014} identify the presence of SNHL in the Wilson-Cowan model and demonstrated that, in a system of strongly coupled components (two interconnected Wilson-Cowan models), it can serve as a mechanism to induce chaos (aperiodic oscillations). Similar implications for chaos induction via SNHL are discussed in the work by Shil'nikov \textit{et al.} \cite{shilnikov_bifurcation_1995} where a simple model of atmospheric circulation is analyzed. A study by Krauskopf \textit{et al.} \cite{krauskopf_excitability_2003} also reports this bifurcation as an organizing center of multipulse excitability in the injection laser. Later codimension-two global bifurcations were studied to understand the transitions to excitability, identifying both homoclinic and heteroclinic structures \cite{wieczorek_bifurcations_2005}. All these studies showcase the birth and disappearance of oscillations via global separatrix-loop bifurcations, such as SNIC and homoclinic in contrast to the more common case of a local Hopf bifurcation. However, heteroclinic bifurcations involving at least one non-hyperbolic equilibrium have received relatively little attention in the literature, particularly in the context of oscillatory dynamics.

 The review \cite[Section 5.2.3]{homburg_chapter_2010} includes a discussion of heteroclinic orbits featuring non-hyperbolic equilibria. The authors remark that a cycle including one saddle-node and one hyperbolic saddle (that we call here `SNICeroclinic') will be generic for two-parameter families, and that a cycle between two saddle-nodes may even perturb to infinitely many periodic orbits. In \cite{grozovski_bifurcations_1996} unfolding of central heteroclinic loop between a saddle node and a saddle (here central SNICeroclinic) is presented, where the phenomenon is referred to as `half-apple' bifurcation. In another paper, Dumortier et al. \cite{dumortier_elementary_1994} consider the unfolding of the heteroclinic loop between a saddle-node and a saddle in two- and three-parameter family featuring saddle-node to saddle heteroclinic loop with non-central connection. This paper focuses on the cyclicity in the system, rather than heteroclinic structures as a route to oscillatory dynamics and does not highlight the applicability of the saddle-node to saddle heteroclinic loop with non-central connection.

In this paper, we are interested in revealing the types and bifurcations of saddle-node separatrix loops that can occur in minimal dimensions, regardless of specific details, helping to classify and understand these transitions. To this end we study the unfolding of a non-central SNICeroclinic loop between one hyperbolic and one non-hyperbolic saddle equilibria, where the non-hyperbolic equilibrium undergoes a saddle-node bifurcation. We carry out our analysis in a two-dimensional setting as this is the simplest scenario in which this bifurcation could occur but note additional features that will appear in higher dimensions. This codimension-three generic bifurcation includes the SNICeroclinic loop as a codimension two bifurcation in its neighborhood.

The analysis below is motivated by two illustrative examples of non-central SNICeroclinic loops arising in dynamical systems (see \autoref{sec:examples}). We employ the method of Shil'nikov variables in our analysis \cite{shilnikov_methods_2001}, reducing the study of a continuous system to studying the associated discrete maps, referred to as a Poincar\'{e} map \cite{poincare_three-body_2017}. In \autoref{sec:ncSNICeroclinic} we present the method of the construction of the Poincar\'{e} maps for the non-central SNICeroclinic loop, splitting the heteroclinic loop into two local and two global dynamics maps. Then in \autoref{sec:unfolding} we consider the dynamics of the Poincar\'{e} map before, after and at the saddle-node bifurcation of the non-hyperbolic equilibrium, combining the findings to reveal the unfolding.  

\section{Examples of non-central SNICeroclinic loops}\label{sec:examples}

A non-central SNICeroclinic loop is a heteroclinic loop connecting a hyperbolic and a non-hyperbolic saddle equilibria, where the orbit approaching the non-hyperbolic saddle is non-central. In this section we give two simple motivating examples of non-central SNICeroclinic loops in two-dimensional dynamical systems. We note there are different types of the phenomenon depending on whether the non-central eigenvalue $\rho$ at the saddle-node is stable or unstable and whether the stable $\lambda_s$ or unstable $\lambda_u$ eigenvalue of the saddle dominates. Table~\ref{tab:cases} highlights four possible cases.
\begin{table}[htbp]
    \centering
    \begin{tabular}{c|cc}
    & $\vert\lambda_s\vert>\vert\lambda_u\vert$ & $\vert\lambda_s\vert<\vert\lambda_u\vert$\\
    \hline
$\rho<0$ &  Type I (stable) & Type II (mixed) \\
$\rho>0$ & Type III (mixed) & Type IV (unstable)
    \end{tabular}
    \caption{Cases of the planar non-central SNICeroclinic loop bifurcation where $\rho$ is the non-central eigenvalue of the saddle-node and $\lambda_s$, $\lambda_u$ are the stable and unstable eigenvalues of the saddle, respectively. In this paper we concentrate on the stable case Type I and show an example of the unstable case (Type IV) which corresponds to time-reversal of Type I.}
    \label{tab:cases}
\end{table}

\subsection{A stable non-central SNICeroclinic loop in a polynomial singular fast-slow system}

\begin{figure}[h!]
    \centering 
    \includegraphics[width=0.8\textwidth]{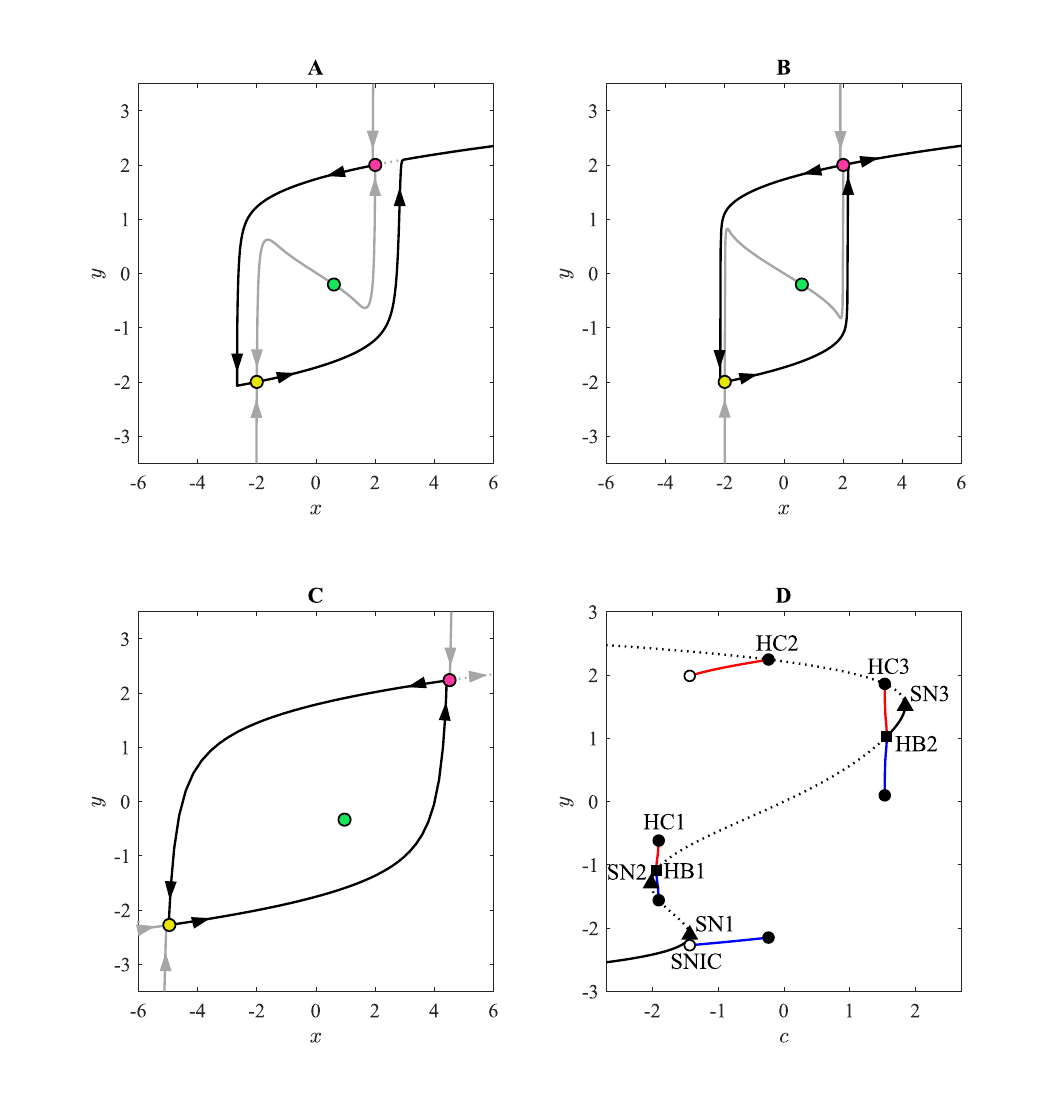}
    \caption{(\textbf{A,B}) Manifolds of the saddle point (pink) and saddle node (yellow) in the polynomial system ({\autoref{eq:poly}) as it approaches the singular limit $(a,b,c)=(0.222\dot{2},1,-0.888\dot{8})$  with $\epsilon=0.1$ (\textbf{A}) and $\epsilon=0.01$ (\textbf{B}). Green dot identifies the location of the unstable node. As $\epsilon \rightarrow 0$ the system approaches the existence of the SNICeroclinic loop. Solid and dotted gray lines identify stable and unstable manifolds, respectively. (\textbf{C}) Non-central SNICeroclinic loop in the polynomial system \autoref{eq:poly} with $\epsilon=1$, $a=0.042$, $b=0.49575$ and $c=-0.85$. Light yellow and pink circles identify saddle-node and saddle, respectively, while green signifies unstable node}. Black line shows the SNICeroclinic loop between the saddle-node and saddle, solid and dotted gray lines identify stable and unstable manifolds, respectively. (\textbf{D}) Continuation in the parameter $c$ in the polynomial system ($\epsilon=1$, $a=0.042$, $b=0.35$), where solid and dotted black lines identify stable and unstable equilibria, respectively; red and blue lines show the maximum and the minimum of the periodic solutions, respectively; black triangle (SN1, SN2, SN3): saddle node; black square (HB1, HB2): Hopf; black circle (HC1, HC2, HC3): homoclinic; white circle (SNIC): saddle node on invariant circle.}
    \label{fig:sniceroclinic1}
\end{figure}

A simple example of a non-central SNICeroclinic loop can be realized in the following two-dimensional polynomial fast-slow system of differential equations: 
\begin{equation}
\begin{aligned}
     \dot x&=\epsilon\big(g(x)-y\big)\\
     \dot y&=\big(x-f(y)\big),
\end{aligned}
\label{eq:poly}
\end{equation}
where $f(x)$ and $g(y)$ are
\begin{align}
    f(y) &=y^3-3y \\
    g(x) &= ax^2+bx+c. 
\end{align}
 The parameter $\epsilon$ represents the separation of timescales between the 'fast' ($y$) and 'slow' ($x$) variables. By suitable choice of $a$, $b$ and $c$, one can ensure that, in the singular limit $\epsilon\rightarrow 0$, there exists a saddle-node $pq_1$ and a saddle $p_2$ equilibria as shown in \autoref{fig:sniceroclinic1}. In this configuration, a separatrix entering $pq_1$ approaches along a non-central direction, and there is also a connection from $pq_1$ to $p_2$.
 One can analytically determine parameters in the singular limit ($\epsilon\rightarrow0$) under which the system exhibits a non-central SNICeroclinic loop. We do so by considering the nullcline geometry. Saddle-node and saddle equilibria must occur at the drop points corresponding to the folds of the critical manifold $x=f(y)$, hence the saddle-node $pq_1$ is at $(x,y)=(-2,-2)$, while the saddle point $p_2$ is at $(x,y)=(2,2)$. For the connections $pq_1 \rightarrow p_2$ and $p_2 \rightarrow p_1$ to exist the $x$-nullcline must be zero at $pq_1$ and $p_2$, respectively. This gives us two equations $-2=4a+3b+c$ and $2=4a+3b+c$ for the $x$-nullcline. For $pq_1$ to be a saddle-node we need the $x$- and $y$-nullclines to be tangent at $pq_1$ meaning that $2ax+b\vert_{pq_1}=-4a+b=1/9$. Solving these three constraints we find a non-central SNICeroclinic loop at $(a,b,c)=(0.222\dot{2},1,-0.888\dot{8})$ for the polynomial system in the singular limit. We have plotted the manifolds of the saddle and saddle nodes in the system at the singular limit for $\epsilon=0.1$ and $\epsilon=0.01$ (see \autoref{fig:sniceroclinic1}(\textbf{A},\textbf{B})). We observe that as $\epsilon \rightarrow 0$ the trajectory that enters the neighborhood of the saddle node (yellow dot) becomes more non-central, whereas the trajectory entering the neighborhood of the saddle (pink dot) aligns more closely with its stable manifold. Consequently, a non-central SNICeroclinic arises in the system for $\epsilon\ll 0$.

The equilibria and connections will persist under small perturbations of $\epsilon$ if we change $(a,b,c)$, meaning there is a curve of parameters $(a(\epsilon),b(\epsilon),c(\epsilon))$, such that a non-central SNICeroclinic loop exists in the system for small enough $\epsilon>0$. Here we demonstrate the existence of the phenomenon at $\epsilon=1$ (see \autoref{fig:sniceroclinic1}(\textbf{C})), which is an example of Type I non-central SNICeroclinic loop. We compute one-parameter bifurcation diagram in $c$, which illustrates a limit cycle born through a SNIC bifurcation, while its destruction occurs via a homoclinic (see \autoref{fig:sniceroclinic1}(\textbf{D})). This represents a novel route to oscillations, where periodic solutions appear and disappear via global bifurcations, namely SNIC and homoclinic.

\subsection{An unstable non-central SNICeroclinic loop in a GTPase activation model}\label{sec:gtp}

If the saddle-node possesses a positive, non-zero eigenvalue, this implies the existence of a linearly unstable direction. We demonstrate this case of an unstable SNICeroclinic loop using the GTPase-tension model of \cite{zmurchok_coupling_2018,zmurchok_local_2023} given by the following ODEs:
\begin{equation}    
\begin{aligned}
    \frac{d L}{d t}&=-\varepsilon\left(L-L_0(G,\phi_1)\right)\\
    \frac{d G}{d t}&=\left(b+f(L)+\gamma \frac{G^n}{1+G^n}\right)\left(G_T-G\right)-G,
\end{aligned}
\label{eq:GTPase}
\end{equation}
where $G$ and $L$ depict the GTPase concentration and cell size, respectively. The function $f(L)$ is a mechanical feedback mechanism from cell deformation to GTPase activity. In \cite{zmurchok_coupling_2018} different forms of $f$ under different assumptions have been considered. For the purposes of our analysis we have selected the Hill function:
\begin{equation}
    f(L)=\beta \frac{L^m}{L_0(G,\phi_2)^m+L^m}. 
\end{equation}
The term $L_0$ represents the rest length that is assumed to decrease from the original length $\ell_0$ depending on the amount of active GTPase in the cell:
\begin{align}
    L_0=L_0(G,\phi)=\ell_0-\phi \frac{G^p}{G_h^p+G^p},
\end{align}
where the parameter $\phi$ is the scaled rate of the GTPase activation. In our analysis we assume that the GTPase concentration and cell size integrate the active GTPase differently, hence we use $\phi_1$ and $\phi_2$ to represent it. The detailed description of the model can be found in \cite{zmurchok_coupling_2018,zmurchok_local_2023}. The parameters used to achieve a non-central SNICeroclinic in the GTPase activation model are in \autoref{tab:gtp_par}.
\begin{table}[]
    \centering
    \begin{tabular}{c|cc}
   Parameter & Value & Description\\
    \hline
$\beta$ &  0.0052 & strength of feedback from tension to GTPase activation\\
$b$ & 0.2530 &  basal activation rate \\
$\gamma$ & 1.6 &  scaled rate of feedback activation \\
$G_T$ & 2 &  mechanical activation constant \\
$\ell_0$ & 1 &  rest length \\
$\phi_1$ & 0.9 & Hill function amplitude  \\
$\phi_2$ & 2 &  Hill function amplitude \\
$G_h$ & 0.4 & half-maximum GTPase activity  \\
$\epsilon$ & 0.1 & rate of contraction  \\
$n$, $p$, $m$ & 4 & Hill coefficients  \\

    \end{tabular}
    \caption{Model parameters for the GTPase activation model \autoref{eq:GTPase}.}
    \label{tab:gtp_par}
\end{table}

\begin{figure}[!h] 
    \centering 
    \includegraphics[width=0.65\textwidth]{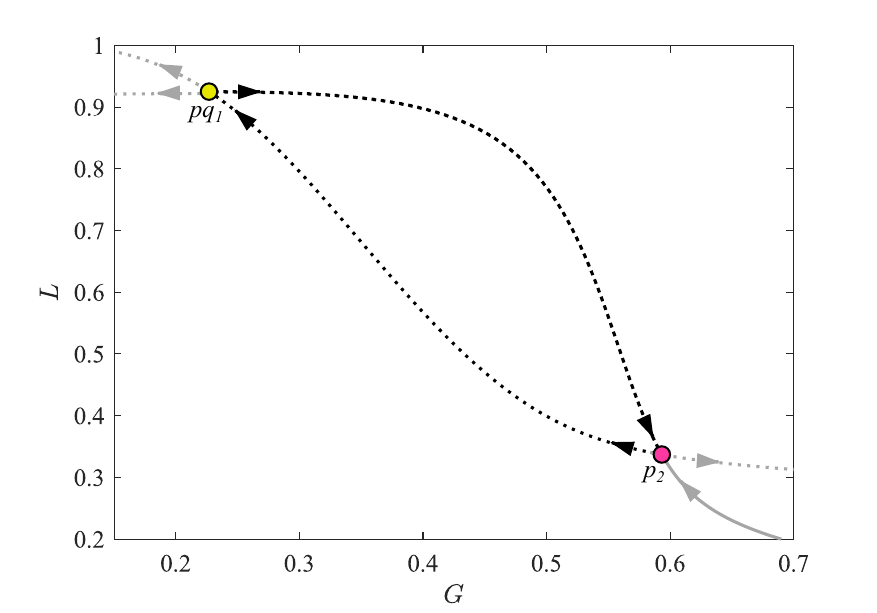} 
    \caption{A non-central SNICeroclinic loop between a saddle-node $pq_1$ with linearly unstable direction and saddle $p_2$ in the GTPase activation model \autoref{eq:GTPase}. Yellow and pink circles identify saddle-node and saddle, respectively. Black dotted line shows the unstable SNICeroclinic loop between the saddle-node and the saddle, solid and dotted gray lines signify the stable and unstable manifolds, respectively.}
    \label{fig:het1}
\end{figure}
The nonhyperbolic point has a positive (unstable non-zero) eigenvalue, and, as a result, the heteroclinic connections will be formed on the intersection of the unstable manifold of the saddle and stable central manifold of the saddle-node as well as the intersection of the stable manifold of the saddle and unstable manifold of the saddle-node (\autoref{fig:het1}). Moreover, the unstable eigenvalue dominates the stable eigenvalue of the saddle. This means we have an unstable (Type IV in Table~\ref{tab:cases}) non-central SNICeroclinic loop. The non-central entrance in the neighborhood of the saddle-node equilibrium arises inherently. Otherwise, the flow will leave the neighborhood of the SNICeroclinic loop, prompted by the unstable separatrix. The below unfolding of a codimension three saddle–node to saddle separatrix-loop bifurcation can be applied to this scenario by reversing time.

\section{The non-central SNICeroclinic loop bifurcation} %
\label{sec:ncSNICeroclinic}

We consider this as an identification problem in general dimensions before giving an unfolding for the planar case.

\subsection{The general non-central SNICeroclinic loop}

Consider an $n$-dimensional dynamical system given by differential equations: 
\begin{equation}
    \frac{du}{dt}=f(u) \;\;\;\;\;\ u \in \mathbb{R}^n , \label{eq:systemRn}  
\end{equation}  
where $f: \mathbb{R}^n \rightarrow \mathbb{R}^n$ is $C^r$ ($r\geq 2$). Assume that \autoref{eq:systemRn} has a nonhyperbolic equilibrium $pq_1$, a hyperbolic saddle $p_2$ and connections between them as in \autoref{fig:schematics}(\textbf{A}).

More precisely, assume $pq_1$ and $p_2$ are connected via a heteroclinic loop formed by the separatrices $\Gamma_1$ and $\Gamma_2$:
\begin{align}
    \Gamma_1 :=& \{z=u(t), t \in \mathbb{R} \;\vert \; u(+\infty)=p_2 \text{ and }  u(-\infty)=pq_1 \}\\
    \Gamma_2 :=& \{z=u(t), t \in \mathbb{R}\; \vert\; u(+\infty)=pq_1 \text{ and }  u(-\infty)=p_2\}. 
\end{align}
We consider the heteroclinic contour formed by  $\Gamma=\overline{\Gamma_1 \cup \Gamma_2}$ and are interested in the invariant sets and bifurcations that appear within some neighborhood $\mathcal{U}$ of $\Gamma$. We say $\Gamma$ is a {\em SNICeroclinic loop} if:
\begin{itemize}
\item[{\bf H1*}]  The system \autoref{eq:systemRn} has a generic saddle-node equilibrium $pq_1$ whose Jacobian $Df(pq_1)$ has a simple eigenvalue $0$ and all other eigenvalues lie on one side of the imaginary axis. The eigenvalues are in general position with no low-order resonances.
\item[{\bf H2*}]  The system has a saddle equilibrium $p_2$ whose Jacobian $Df(p_2)$ has leading eigenvalues $\lambda_u$ and $\lambda_s$, with $\Re(\lambda_u)>0$ and $\Re(\lambda_s)<0$. The eigenvalues are in general position with no low-order resonances.
\item[{\bf H3*}]
There is a single trajectory $\Gamma_1$ from $pq_1$ to $p_2$ that is a transversal intersection of the center-unstable manifold of $pq_1$ and the stable manifold of $p_2$.
\item[{\bf H4*}] 
There is a single trajectory $\Gamma_2$ from $p_2$ to $pq_1$ that is a transversal intersection of the unstable manifold of $p_2$ and the center-stable manifold of $pq_1$.
\end{itemize}
The SNICeroclinic loop is \textit{non-central} if:

\begin{itemize}
\item[{\bf H5*}]
$\Gamma_2$ approaches $pq_1$, or $\Gamma_1$ leaves $pq_1$ in a non-central direction. 
\end{itemize}
 We simplify these assumptions and consider an unfolding in the special case $n=2$ in the next section. Challenges associated with the general unfolding problem are discussed in \autoref{sec:discussion}.

\subsection{Unfolding the planar case}

Now consider the parametrized planar dynamical system: 
\begin{equation}
    \frac{du}{dt}=f(u,\mu) \;\;\;\;\;\ (u,\mu) \in \mathbb{R}^2 \times \mathbb{R}^3, \label{eq:system}  
\end{equation}  
where $f: \mathbb{R}^2\times\mathbb{R}^3 \rightarrow \mathbb{R}^2$ is $C^r$ ($r\geq 2$).
We will unfold this using the Kuznetsov Normal Form Theorem \cite{kuznetsov_elements_2023} to locally approximate the system near the non-central SNICeroclinic loop and then proceed to construct a Poincar\'{e} map. Construction of the Poincar\'{e} map involves a reduction to a co-dimension one surface where we preserve local dynamics of the original system \cite{broer_chapter_2010}. We follow the method from the work of Dumortier \textit{et al.} \cite{dumortier_elementary_1994} but we aim to make the assumptions and calculations more explicit and generalizable to higher-dimensional cases.
\begin{figure}[htbp] 
    \centering 
    \includegraphics[width=0.7\textwidth]{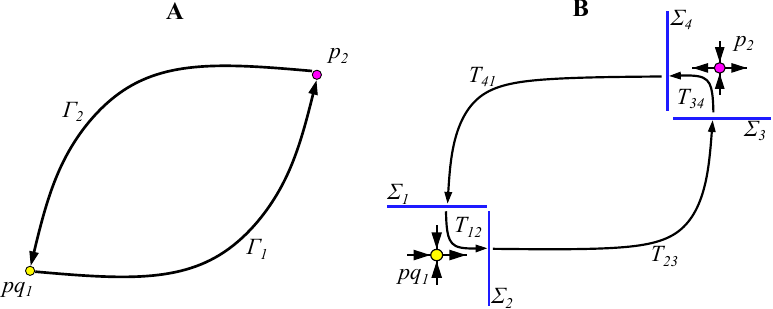} 
    \caption{\textbf{(A)} Non-central SNICeroclinic loop $\Gamma=\overline{\Gamma_1 \cup \Gamma_2}$ between saddle-node $pq_1$ and saddle $p_2$. \textbf{(B)} Overview of the construction of the Poincar\'{e} map for the planar case in the vicinity of the heteroclinic loop $\Gamma$  between nonhyperbolic equilibrium ($pq_1$, yellow) and hyperbolic saddle ($p_2$, pink) with established cross-sections ($\Sigma_1$, $\Sigma_2$, $\Sigma_3$, $\Sigma_4$) and corresponding ``connection'' maps $T_{12}$, $T_{23}$, $T_{34}$, $T_{41}$.}
    \label{fig:schematics}
\end{figure}
In the planar case $n=2$ we assume there is a $\Gamma$ with neighborhood $\mathcal{U}$ such that {\bf H1*-H4*} become the following hypotheses for $\mu=0$:
\begin{itemize}
\item[{\bf H1}]  The system \autoref{eq:system} has an equilibria $pq_1$ whose Jacobian $Df(pq_1)$  has real eigenvalues $\rho\neq0$ and $0$. This is a generic saddle-node of codimension 1.
\item[{\bf H2}]  The system has an equilibrium $p_2$ whose Jacobian $Df(p_2)$ has two simple real eigenvalues: $\lambda_u>0$ and $\lambda_s<0$, such that $\vert\lambda_s \vert\neq \vert\lambda_u \vert$ and there are no low order resonances. 
\item[{\bf H3}]
The one-dimensional center $W^c(pq_1)$ contains a trajectory $\Gamma_1$ that limits to $p_2$.
\item[{\bf H4}] 
The  one-dimensional unstable $W^u(p_2)$ contains a trajectory $\Gamma_2$ that limits to $pq_1$.

\end{itemize}

The assumptions \textbf{H1-H4} define the existence of the \textit{SNICeroclinic loop} in the 2-\\dimensional system. While to ensure that the loop is \textit{non-central} we specify: 

\begin{itemize}
    \item[{\bf H5}] 
$\Gamma_2$ approaches $pq_1$ in a non-central direction. 
\end{itemize}
If \textbf{H5} is not satisfied this will be generic in two-parameter families \cite{grozovski_bifurcations_1996}.In addition we assume 
\begin{itemize}
\item[\bf H6] The parametrization by $\mu$ is in general position and it unfolds the bifurcation at $\mu=0$.
\end{itemize}

The assumption \textbf{H6} ensures that we have enough parameters and they are in correct position to unfold. 

Even in the planar case, there are several types of non-central SNICeroclinic loop in the plane, depending on the sign of $\rho$ and $\lambda_s, \lambda_u$, as outlined in Table~\ref{tab:cases}. In this paper, we concentrate on stable case Type I which results in bifurcation to attracting periodic orbits. The unstable case Type IV corresponds to the case with saddle-node with linearly unstable direction, the existence of which was demonstrated in the GTPase activation model (\autoref{fig:het1}). The unfolding of the mixed cases Type II and Type III can be found in \cite{dumortier_elementary_1994} and includes both stable and unstable periodic orbits that meet at a saddle node of periodic orbits.

We perform geometric decomposition of the heteroclinic loop using local and global maps to derive the conditions for different bifurcations to occur varying the parameters $\mu$. We argue that the considered non-central SNICeroclinic loop is a codimension-three structure, where one parameter is responsible for the saddle-node bifurcation of $pq_1$, and two other parameters control the bifurcation of the trajectories $\Gamma_1$ and $\Gamma_2$. 
Let $\mu_1$ be the unfolding parameter of the saddle-node bifurcation of $pq_1$ and let $\mu_2$ and $\mu_3$ be the splitting parameters for separatrices $\Gamma_1$ and $\Gamma_2$, respectively. We separate the heteroclinic loop $\Gamma$ into 4 segments with distinct qualitative dynamics (see \autoref{fig:schematics}(\textbf{B})). The maps $T_{12}$ and $T_{34}$ are local maps around points $pq_1$ and $p_2$, respectively, while the maps $T_{23}$ and $T_{41}$ are the global transition maps from neighborhood $U_1$ to $U_2$ and from $U_2$ to $U_1$, respectively.

 \subsection{Local map around the saddle-node $pq_1$: $T_{12}$}
 
According to the center manifold theory there exists a center manifold that governs the flow on a small neighborhood of $pq_1$ \cite{carr_applications_1981}. We use the saddle-node normal form to define the unfolding of local dynamics around the box neighborhood $U_1:= \{ (x,y) \in \mathbb{R}^2 \;\vert\; \vert x \vert, \vert y \vert \leq \delta \}\subset \mathcal{U}$ of the critical point $pq_1$, where $\delta \ll 1$ \cite{ermentrout_parabolic_1986,baesens_interaction_2013}. Under this choice of variables, the flow is as follows: 
\begin{align}
    \dot x =& x^2-\mu_1 \label{T4x}\\
    \dot y =& \rho y \label{T4y},
\end{align}
where $\mu_1$ is the unfolding parameter, and $\rho$ is the non-central eigenvalue. We note that using the above change of variables we have topological equivalence between the normal form and the original system. The topologically equivalent normal form is homeomorphic to the original system and maintains the time orientation \cite{kuznetsov_elements_2023}, which is sufficient for the analysis here. To consider the degree of asymmetry of the two branches we direct the reader to Glendinning and Simpson  \cite{glendinning_normal_2022} as well as  Ilyashenko \textit{et al.} \cite{ilyashenko_finitely-smooth_1991}, where the normal form under smooth equivalence has a cubic term. Here we assume that the cubic term is zero, as with the cubic term, the two equilibria that bifurcate at the saddle-node bifurcation no longer have symmetrically related eigenvalues with opposite signs, but we note that this does not affect the branching or stability of the invariant sets that we consider here, which are already determined by the topological normal form.

We use the method of Shil'nikov variables \cite{shilnikov_methods_1998} to find the transition map from $\Sigma_1$ to $\Sigma_2$. Let us define the cross-sections as follows:
\begin{align}
        \Sigma_1 :=& \{ (x,y) \in \mathbb{R}^2\; \vert\; |x|<\delta, y = \delta \}\\
        \Sigma_2 :=& \{ (x,y) \in \mathbb{R}^2\; \vert\; |y|<\delta, x = \delta \}. 
\end{align}
We suppose that at time $t=0$ the flow intersects $\Sigma_1$ at $x=x_0>0$ and at time $t=\tau$ it intersects $\Sigma_0$ at $y=y_1>0$. We have three cases:$\mu_1>0$, $\mu_1<0$ and $\mu_1=0$. 
\begin{figure}[h] 
    \centering 
    \includegraphics[width=1\textwidth]{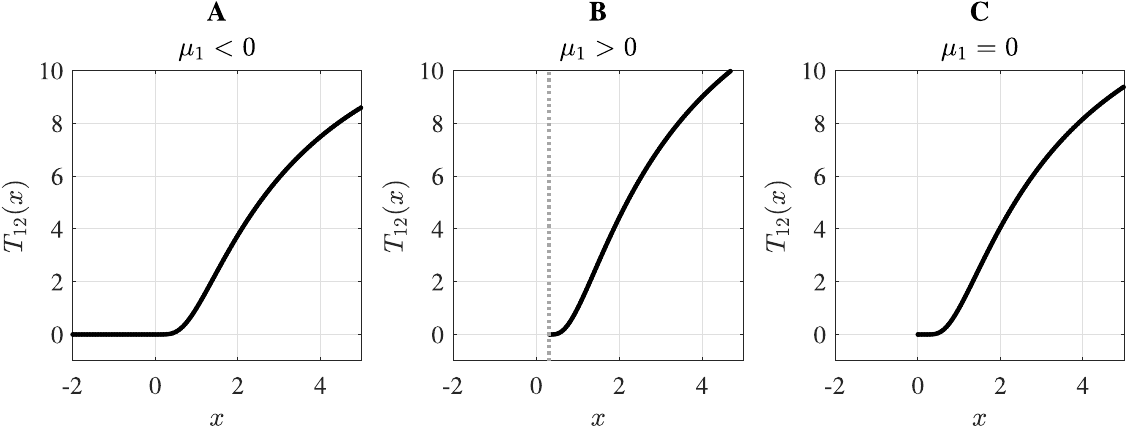} 
    \caption{Representation of the connecting maps  $T_{12}$ for (\textbf{A}) $\mu_1<0$, (\textbf{B}) $\mu_1>0$, (\textbf{C}) $\mu_1=0$.} 
    \label{fig:conmapst12}
\end{figure}
\begin{itemize}
    \item \textit{ Case 1: $\mu_1<0$}
    
    For $\mu_1<0$ the system of equations \eqref{T4x} and \eqref{T4y} has no real fixed points, i.e. the saddle-node disappears for $\mu_1<0$. 
    The general solution of the system of the equations \eqref{T4x} and \eqref{T4y} in this case is the following: 
\begin{align}
    t+C_1=&\frac{1}{\sqrt{|\mu_1}|}\arctan({\frac{x}{\sqrt{|\mu_1|}}})\\
    y=&C_2e^{t\rho},
\end{align}
where $C_1$ and $C_2$ are constants. Using initial conditions $(x_0, \delta)\in \Sigma_1$ we obtain 
\begin{align}
    C_1=&\frac{1}{\sqrt{|\mu_1|}}\arctan(\frac{x_0}{\sqrt{|\mu_1|}})\\
    C_2=&\delta.
\end{align}
At $t=\tau$ we assume $(x,y)\in \Sigma_2$, thus we have $x=\delta$ and $y=y_1$, where $|y_1|<\delta$, yielding: 
\begin{align}
    \tau=& \frac{1}{\sqrt{|\mu_1|}}\left[\arctan({\frac{\delta}{\sqrt{|\mu_1|}}})-\arctan(\frac{x_0}{\sqrt{|\mu_1|}})\right] \label{T4x0}\\
    y_1=&\delta e^{\tau \rho}. \label{T4y1}
\end{align}
As $\tau$ must be positive we require $x_0 = \{x\in\mathbb{R} \vert x<\delta\}$. Substituting the expression for $\tau$ in equation \ref{T4x0} into \ref{T4y1} we obtain the expression for transition map $T_{12}$ (see \autoref{fig:conmapst12}(\textbf{A})): 
\begin{equation}
    T_{12}= \delta \exp{\left[\frac{\rho}{\sqrt{|\mu_1|}}(\arctan({\frac{\delta}{\sqrt{|\mu_1|}}})-\arctan(\frac{x_0}{\sqrt{|\mu_1|}}))\right]}
\end{equation}
At small values of $x<0$ the map $T_{12}$ approaches zero, but never equal to it.

\item \textit{ Case 2: $\mu_1>0$}

    The linearized system for $\mu_1>0$ has two fixed points, namely a stable $(-\sqrt{\mu_1}, 0)$ and a hyperbolic saddle $(\sqrt{\mu_1}, 0)$. We refer to the saddle point as $p_1$. 
    The general solution of the system of the equations \eqref{T4x} and \eqref{T4y} with positive $\mu_1$ is 
\begin{align}
    t+C_1=&\frac{1}{2\sqrt{\mu_1}}\ln \left|{\frac{x-\sqrt{\mu_1}}{x+\sqrt{\mu_1}}}\right|\\
    y=&C_2e^{t\rho},
\end{align}
where $C_1$ and $C_2$ are constants. Using initial conditions $(x_0, \delta)\in \Sigma_1$ we obtain 
\begin{align}
    C_1=&\frac{1}{2\sqrt{\mu_1}}\ln \left|{\frac{x_0-\sqrt{\mu_1}}{x_0+\sqrt{\mu_1}}}\right|\\
    C_2=&\delta.
\end{align}
At $t=\tau$ we hit $(\delta,y_1)\in \Sigma_2$, which yields: 
\begin{align}
    \tau=& \frac{1}{2\sqrt{\mu_1}}\left[\ln \left|{\frac{\delta-\sqrt{\mu_1}}{\delta+\sqrt{\mu_1}}}\right|-\ln \left|{\frac{x_0-\sqrt{\mu_1}}{x_0+\sqrt{\mu_1}}}\right|\right] \label{T4x02}\\
    y_1=&\delta e^{\tau \rho}. \label{T4y12}
\end{align}
As $\tau$ must be positive we require $x_0 = \{x\in\mathbb{R} \vert \sqrt{\mu_1}\leq x<\delta \}$, putting restrictions on $\Sigma_1$. Substituting the expression for $\tau$ in equation \eqref{T4x02} into \eqref{T4y12} we obtain the expression for transition map $T_{12}$ for $\mu_1>0$ (\autoref{fig:conmapst12}(\textbf{B})): 
\begin{equation}
    T_{12}= \delta \exp\left[{\frac{1}{2\sqrt{\mu_1}}\left[\ln \left|{\frac{\delta-\sqrt{\mu_1}}{\delta+\sqrt{\mu_1}}}\right|-\ln \left|{\frac{x-\sqrt{\mu_1}}{x+\sqrt{\mu_1}}}\right|\right]}\right].
\end{equation}
Note, we restrict the domain of the 'connecting' map $T_{12}$ for $\mu_1>0$ to be $x \in (\sqrt{\mu_1,\delta})$.

\item \textit{Case 3: $\mu_1=0$} 

When $\mu_1=0$, the linearized system has one non-hyperbolic saddle point $(0,0)$. 

The general solution of the system of the equations \ref{T4x} and \ref{T4y}: 
\begin{align}
    t+C_1=&-\frac{1}{x}\\
    y=&C_2e^{t\rho},
\end{align}
where $C_1$ and $C_2$ are constants. Using initial conditions we obtain 
\begin{align}
    C_1=&-\frac{1}{x_0}\\
    C_2=&\delta.
\end{align}
At $t=\tau$ we have $x=\delta$ and $y=y_1$, yielding: 
\begin{align}
    \tau=& \frac{1}{x_0}-\frac{1}{\delta}\label{T4x03}\\
    y_1=&\delta e^{\tau \rho}. \label{T4y13}
\end{align}
As $\tau$ must be positive we require $x_0 = \{x\in\mathbb{R} \vert 0\leq x<\delta \}$. This imposes restriction to the $\Sigma_1$. Substituting the expression for $\tau$ in equation \ref{T4x03} into \ref{T4y13} we obtain the expression for transition map $T_{12}$ (see \autoref{fig:conmapst12}(\textbf{C})): 
\begin{equation}
    T_{12}(x)= \delta \exp\left[{\rho\left(\frac{1}{x}-\frac{1}{\delta}\right)}\right]
\end{equation}
\end{itemize}

\subsection{Local map around the saddle point $p_2$: $T_{34}$}

Using the Sternberg linearization theorem \cite{stowe_linearization_1986}, we assume there is a $C^2$ coordinate change that linearizes the vector field and hence flattens the stable and unstable manifolds within a small box neighborhood $U_2$ of $p_2$.

\begin{figure}[h] 
    \centering 
    \includegraphics[width=1\textwidth]{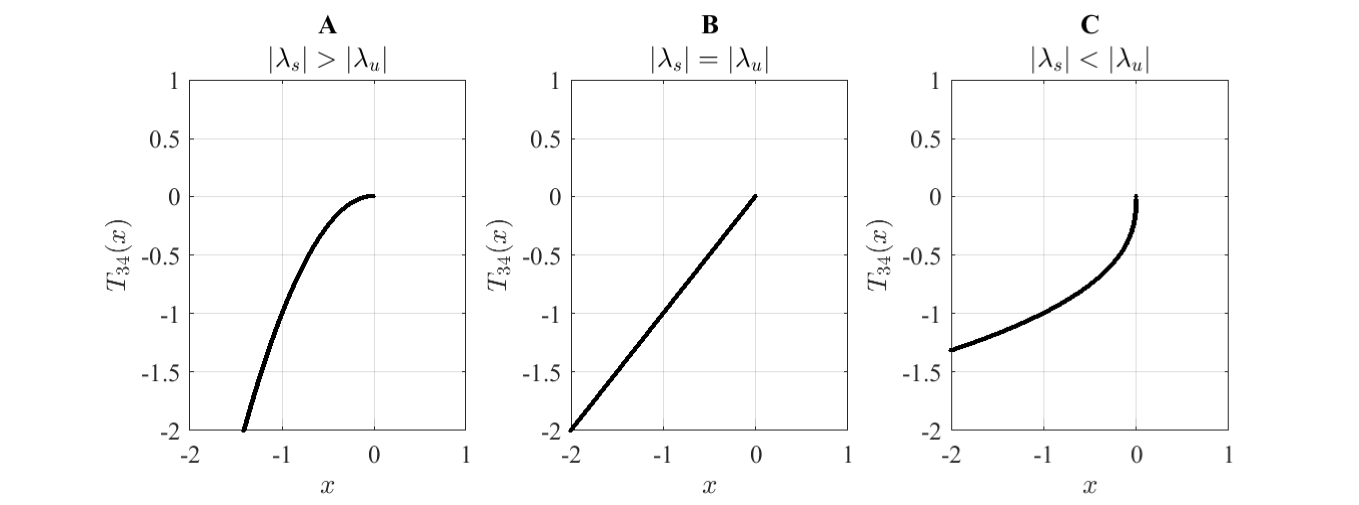} 
    \caption{Representation of the connecting maps  $T_{34}$ for \textbf{(A)} $|\lambda_s|>|\lambda_u|$, \textbf{(B)} $|\lambda_s|=|\lambda_u|$, \textbf{(C)} $|\lambda_s|<|\lambda_u|$.} 
    \label{fig:conmapst34}
\end{figure}

We set local coordinate system $(x,y)^T$ within $U_2$, such that the linear system can be written: 
\begin{align}
    \dot x =& \lambda_u x  \label{xT2}\\
    \dot y =& \lambda_s y.  \label{yT2}
\end{align}
In this coordinate system, the dynamics in the local stable $W^s_{loc}(p_2)$ and unstable manifolds $W^u_{loc}(p_2)$ correspond to the uncoupled subsystems associated with eigenvalues of $\lambda_s$ and $\lambda_u$, respectively. 

Under the construction, $x$ ($y$) is in the direction of unstable (stable) manifolds. Without loss of generality, choose a box $U_2$ in the form: 
\begin{equation}
    U_2 := \{ (x,y) \in \mathbb{R}^2 \;\vert\; \vert x \vert, \vert y \vert \leq \epsilon \},\;\;\;\; \text{for some} \;\; \epsilon \ll 1.
\end{equation}
We can now define the Poincar\'{e} sections on $\partial U_2$:

\begin{align}
        \Sigma_3 :=& \{ (x,y) \in \mathbb{R}^2\; \vert\; |x|<\epsilon,\; y=-\epsilon \}\\
        \Sigma_4 :=& \{ (x,y) \in \mathbb{R}^2\; \vert\; |y|<\epsilon,\; x=-\epsilon \}.
\end{align}

These sections are transverse to the flow provided $|\epsilon|$ is small enough. Clearly trajectories that remain for all positive times in the neighborhood $U_2$ of the heteroclinic loop must intersect $\Sigma_3$ and $\Sigma_4$. Using the method of Shil'nikov variables \cite{shilnikov_methods_2001}, we can find a transition map $T_{34}:\Sigma_3 \rightarrow \Sigma_4$ near a given equilibrium by solving a boundary value problem for the system of equations \eqref{xT2} and \eqref{yT2}. We emphasize these boundary conditions ensure that the flow is to stay in the vicinity of $\Gamma$.

In the neighborhood of $U_2$ let the trajectory start at $(x_0, -\epsilon)\in \Sigma_3$ and end at $(-\epsilon, y_1)\in \Sigma_4$. The general solution to the system of equations \eqref{xT2} and \eqref{yT2} is 
\begin{align}
    x = &C_1 e^{\lambda_u t}\\
    y = &C_2 e^{\lambda_s t},
\end{align}
where $C_1, C_2 \in \mathbb{R}$. Using initial conditions $(x,y)=(x_0,-\epsilon)$ at $t=0$, we obtain $C_1=x_0$ and $C_2=-\epsilon$. Consider the system at $t=\tau$: 
\begin{align}
    -\epsilon=&x_0 e^{\lambda_u \tau} \label{x0T2}\\
    y_1=&-\epsilon e^{\lambda_s \tau}. \label{y1T2}
\end{align}
Rearranging the terms in equation \ref{x0T2} we obtain: 
\begin{align}
    \tau = \frac{1}{\lambda_u}\ln{\frac{\epsilon}{|x_0|}} \label{tauT2}. 
\end{align}
We note that the expression for $\tau$ is defined only for $-\epsilon<x_0\leq0$. 
Substituting the expression for $\tau$ in equation \eqref{tauT2} into equation \eqref{y1T2} we obtain the the transition map $T_{34}$ : 
\begin{equation}
    T_{34}(x) = -\epsilon^{1+\lambda_s/\lambda_u} |x|^{-\lambda_s/\lambda_u},  
\end{equation}
which is defined for $-\epsilon<x\leq0$. We note that the system's evolution is defined by the point at which the trajectory enters the local neighborhood of $p_2$, specifically: 
\begin{itemize}
    \item if the trajectory enters $\Sigma_3$ at $-\epsilon<x_0<0$, it will stay in the neighborhood $\mathcal{U}$ of $\Gamma$ after leaving the local neighborhood of $p_2$
    \item if the trajectory enters $\Sigma_3$ at $x_0=0$, it asymptotically approaches $p_2$.
    \item if the trajectory enters $\Sigma_3$ at $x_0>0$, it will leave the neighborhood $\mathcal{U}$ of $\Gamma$ without passing though the saddle point $p_2$.
\end{itemize}

Furthermore, based on the magnitude of $\lambda_s/\lambda_u$ the shape of the connecting map $T_{34}$ changes (see \autoref{fig:conmapst34}). The magnitude of the eigenvalue of the saddle point tells us how fast the flow travels in that direction \cite{kuznetsov_elements_2023}. In our analysis we only consider the case when $\lambda_s>\lambda_u$, depicted in \autoref{fig:conmapst34}(\textbf{A}).

\subsection{Global return maps}

To define return maps, we need to define maps from $\Sigma_2$ to $\Sigma_3$ and from $\Sigma_4$ to $\Sigma_1$.
We separate the flow into two parts, namely flow that goes from the neighborhood of $pq_1$ and ends at the neighborhood of $p_2$ and flow that starts at the neighborhood of $p_2$ and goes to the neighborhood of $pq_1$ (\autoref{fig:schematics}(\textbf{B})).  To approximate the flow from $\Sigma_2$ to $\Sigma_3$ we use first-order approximation of the path for the system in \eqref{eq:system} at $x_0$ ($y$ component of the cross-section $\Sigma_2$), neglecting higher-order terms: 
\begin{align}
    T_{23}(x)=a_1 x +\mu_2,
\end{align}
 where $a_1$ is a constant and $\mu_2$ is an unfolding parameters. Similarly, the flow from $\Sigma_4$ to $\Sigma_1$ is approximated as 
\begin{align}
    T_{41}(x)=a_2 x +\mu_3,
\end{align}
 where $a_2$ is a constant and $\mu_3$ is an unfolding parameter. As in the neighborhood of $pq_1$ we consider only non-negative values, while in the neighborhood of $p_2$ we restrict the map to negative values, we use global maps $T_{23}$ and $T_{41}$ to preserve the transition between positive and negative values by fixing $a_1=a_2=-1$. 
\subsection{Return maps}
 We ensure that (\textbf{H6}) is satisfied by assuming (a) we have a parameter $\mu_1$ that topologically unfolds the saddle-node $pq_1$ (b) we have a parameter $\mu_2$ that unfolds the connecting orbit from $pq_1$ to $p_2$ by introducing a constant term in the map $T_{23}$ from between $\Sigma2$ and $\Sigma_3$ and (c) we have a parameter $\mu_3$ that unfolds the connecting orbit from $q_2$ to $pq_1$ by introducing a constant term in the map $T_{41}$ from between $\Sigma_4$ and $\Sigma_1$. The maps $T_{23}$ and $T_{41}$ can be assumed affine to lowest order. Using the connecting maps $T_{12}$, $T_{23}$, $T_{34}$, $T_{41}$  we can define return maps to cross-section $\Sigma_1$ and $\Sigma_3$: 
 \begin{align}
     R_1 &=T_{41}\circ T_{34}\circ T_{23}\circ T_{12}(x) \text{ to } \Sigma_1,\\
     R_3 &=T_{23}\circ T_{12}\circ T_{41}\circ T_{34}(x) \text{ to } \Sigma_3.
 \end{align}

\subsection{Analysis method}
 
We observe that in our system with a saddle-node and a saddle point the behavior in the neighborhood $\mathcal{U}$ of $\Gamma$ depends on the behavior of separatrices $\Gamma_1$ and $\Gamma_2$. In \cite{shashkov_complex_1996}, a lemma was proposed that showed that in the unfolding of the heteroclinic loop between two saddle-foci any orbit that stays in the neighborhood of the heteroclinic loop lies in the closure of the union of separatrices connecting the points. Particularly, it implies that the behavior of the trajectories, which start within the neighborhood of heteroclinic loop, is determined by the behavior of the separatrices that are associated with each saddle. We apply this theoretical finding to our analysis. Specifically, we discuss the conditions to induce a specific qualitative behavior in the system, then validate them by considering the stable fixed points of the return maps. The return maps' fixed points allows us to identify whether a separatrix leaves the neighborhood $\mathcal{U}$ of $\Gamma$ or not. To this end we split our analysis of the parameter space of $\mu$ into three parts: $\mu_1=0$, $\mu_1<0$, $\mu_1>0$ -, as $T_{12}$ is split into three cases. For each of the three cases for $\mu_1$ we will analyse the dynamics in the system by considering return maps under variation of splitting parameters $(\mu_2, \mu_3)$. We start by fixing the parameter $\mu_1=0$.

\section{Unfolding the non-central SNICeroclinic}
\label{sec:unfolding}

\subsection{Qualitative behavior changes at the saddle-node bifurcation ($\mu_1=0$)}
 
\begin{figure}[!h] 
    \centering 
    \includegraphics[width=0.99\textwidth]{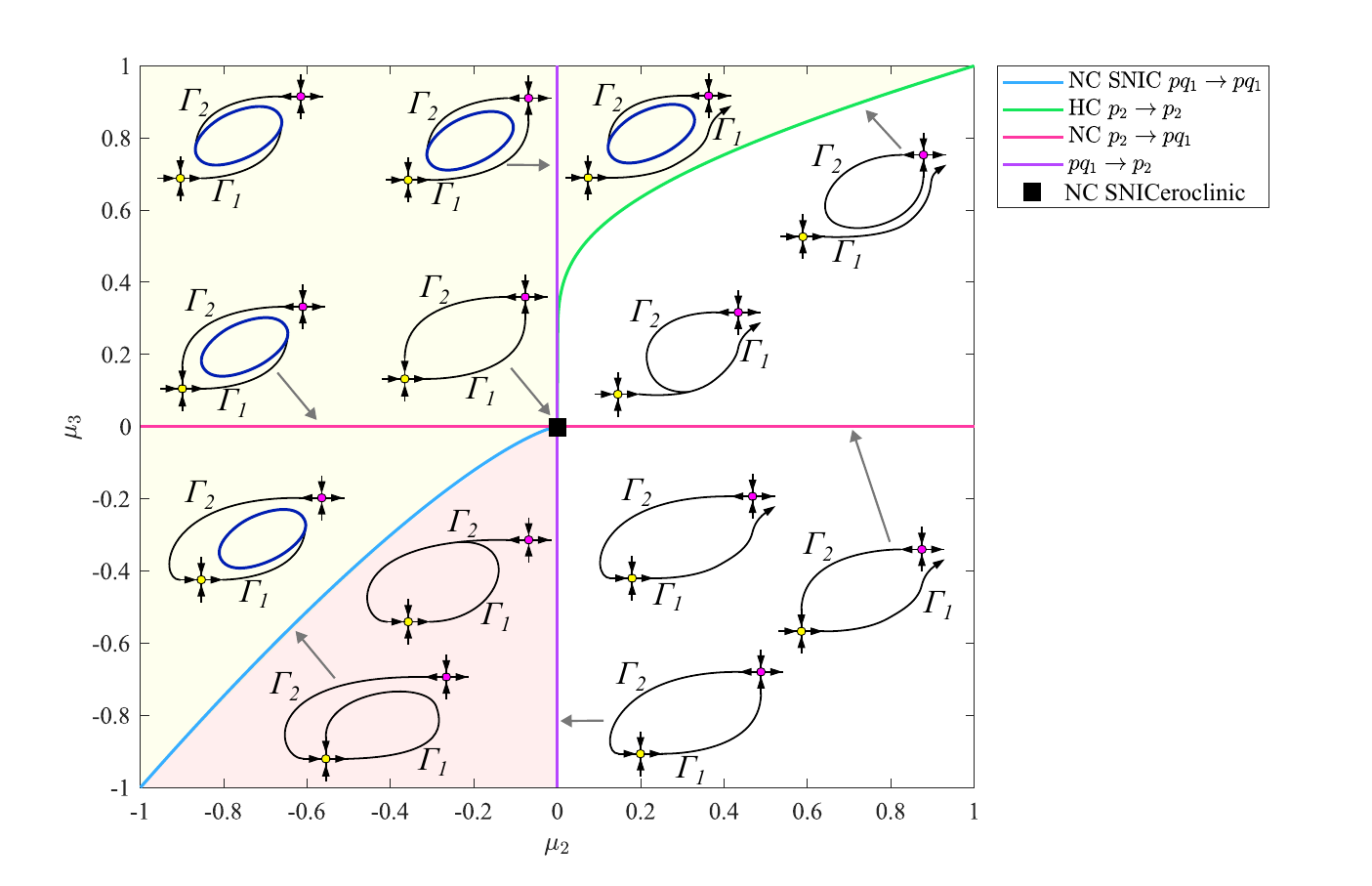} 
    \caption{Slice of the unfolding of the heteroclinic loop between non-hyperbolic and hyperbolic points for fixed $\mu_1=0$, while the parameters $\mu_2$ and $\mu_3$ are changed. Black dot signifies the non-central heteroclinic loop, the curves separate the parameter regions of different dynamics as described in the legend. Light yellow and white colored regions signify the presence and absence of a periodic orbit, respectively, while light red colored region indicates the existence of the central SNIC. The dynamics of the system is also illustrated via phase portraits, depicting the behavior of the separatrices $\Gamma_1$ and $\Gamma_2$. Dark blue closed curves identify the presence of a periodic orbit. SNIC: saddle-node on invariant circle, HC: homoclinic loop, NC HET: non-central heteroclinic loop.} 
    \label{fig:mu1neut}
\end{figure}

When \\$\mu_1=0$ the system has two fixed equilibria - a saddle-node and a saddle. Varying the splitting parameters $\mu_2$ and $\mu_3$ we aim to explore the range of distinct dynamic behaviors. It includes non-central and central heteroclinic connection between $pq_1$ and $p_2$, homoclinic loop with $p_2$, non-central and central SNIC with $pq_1$, existence of periodic solutions with diverse behavior of separatrices $\Gamma_1$ and $\Gamma_2$ as well as destruction of invariant structures. All scenarios on the parameter space $(\mu_2, \mu_3)$ are demonstrated in \autoref{fig:mu1neut}. In the following subsections, we will explore the necessary conditions for maintaining the persistence of a specific dynamic state within the system.
    
\subsubsection{At the non-central SNICeroclinic bifurcation}

When defining the Poincar\'{e} maps we assumed that at $\mu=(0,0,0)$ there exists a heteroclinic loop between $pq_1$ and $p_2$, which we define as the location of the non-central SNICeroclinic bifurcation. In this case each of the four cross-section maps has a stable fixed point at $x=0$ (see \autoref{fig:schematics}), indicating the presence of a bidirectional heteroclinic connection between equilibria following a trajectory entering the saddle-node from the non-central direction. Moreover, given the way we defined the maps if we fix $\mu_2=0$ and vary $\mu_3$ we preserve  heteroclinic connection $pq_1 \rightarrow p_2$ formed by $\Gamma_1$. Depending on the sign of $\mu_3$ we either ensure a persistence of a periodic orbit for $\mu_3>0$ or formation of central heteroclinic connection $p_2 \rightarrow pq_1$. The later case will lead to the persistence of central heteroclinic loop. Similarly, if we fix $\mu_3=0$ and vary $\mu_2$ we can preserve the persistence of non-central heteroclinic connection $p_2 \rightarrow p1$, formed by $\Gamma_2$. Then, if $\mu_2<0$ the separatrix $\Gamma_1$ stays in the neighborhood $\mathcal{U}$, tending to a periodic orbit. If $\mu_2>0$, the separatrix $\Gamma_1$ leaves the neighborhood $\mathcal{U}$.

\subsubsection{Homoclinic loop $p_2 \rightarrow p_2$}\label{sec:hom1}

The homoclinic orbit exists on the intersection of stable and unstable manifolds of the saddle $p_2$. To ensure this condition holds for our return maps, the trajectories in the neighborhood $\mathcal{U}$ of $\Gamma$ are to enter $\Sigma_3$ at a point $(0,\epsilon)$, thus we require the fixed point of the return map $R_3$ to be 0: 
\begin{align}
    0 & = R_3(0) \\
    0 & = a_1T_{12}(\mu_3)+\mu_2\\
    \mu_2 & = -a_1T_{12}(\mu_3),
\end{align}
giving us the conditions on the parameters $\mu_2$ and $\mu_3$. Additionally, the separatrix $\Gamma_2$ must not form a heteroclinic connection (both central and non-central), hence $\mu_3>0$. Under these conditions, the separatrix $\Gamma_2$ forms a homoclinic loop, while $\Gamma_1$ leaves the neighborhood. If we set $\mu_2 > -a_1T_{12}(\mu_3)$ the return map to the cross-section $\Sigma_3$ will lose a stable fixed point, indicating that the separatrix $\Gamma_2$ leaves the neighborhood $\mathcal{U}$, uniting with the separatrix $\Gamma_1$. On the other hand, if $\mu_2 < -a_1T_{12}(\mu_3)$, the return maps have non-zero stable fixed points, demonstrating that the separatrix $\Gamma_2$ tends to a periodic orbit in the neighborhood $\mathcal{U}$.
Hence, the curve $\mu_2 = -a_1T_{12}(\mu_3)$ for the homoclinic loop $p_2 \rightarrow p_2$ is the mechanism that induces transition between persistence of periodic solutions in the system and their destruction. 

\subsubsection{Non-central saddle-node on invariant circle (SNIC)  $pq_1 \rightarrow pq_1$}

The behavior of manifolds in the context of SNIC is a complicated questions, as the relative position of the manifolds varies \cite{homburg_chapter_2010}.
In the case when we have non-central SNIC stable and unstable central manifolds conjoin, which is also consistent with what occurs at the SNHL bifurcation. To observe saddle-node on invariant circle (SNIC) $pq_1 \rightarrow pq_1$ we require the trajectories to enter the neighborhood of $pq_1$ at $\Sigma_1=(0,\delta)$, hence $R_1$ must have a fixed point at $x=0$:
\begin{align}
    0&=R_1(0),\\
    0&=a_2T_{34}(\mu_2)+\mu_3,\\
    \mu_3&= -a_2T_{34}(\mu_2), 
\end{align}
which defines the relationship between the splitting parameters. Additionally, we require the separatrix $\Gamma_1$ not to leave the neighborhood $\mathcal{U}$ of $\Gamma$, hence $\mu_2<0$. This ensures that there exists a connection formed by separatrix $\Gamma_1$ from $pq_1$ to itself. If we set $\mu_3<-a_2T_{34}(\mu_2)$, the non-central connection from $pq_1$ to itself will become central. 
In this case the separatrix $\Gamma_1$ forms a SNIC that enters the neighborhood of $pq_1$ from central direction direction. The other way around $\mu_3>-a_2T_{34}(\mu_2)$ ensures the persistence of the limit cycle in the system and the separatrix $\Gamma_1$ stays in the neighborhood $\mathcal{U}$ tending to the periodic orbit. 

\subsubsection{Birth and death of periodic solutions} 
 In the parameter space ($\mu_2$, $\mu_3$) with $\mu_1=0$, the SNIC, homoclinic, and SNICeroclinic bifurcations together govern the transition from the persistence of periodic orbits to trajectories leaving the neighborhood $\mathcal{U}$ along the separatrix $\Gamma_1$. The case $\mu_3<0$ allows to maintain central heteroclinic connection $p_2 \rightarrow pq_1$, while for $\mu_3>0$ the connection is broken. In the parameter space where the periodic solutions exist the choice of the parameters decides the the qualitative behavior of the separatrices and their interactions with the periodic orbit. Namely, if we are in the second quadrant ($\mu_2<0$ and $\mu_3>0$), both separatrices will tend to the periodic orbit. If we cross the threshold and enter the first quadrant ($0<\mu_2<-a_1T_{12}(\mu_3))$), $\Gamma_2$ continues to tend to the periodic orbit, while  $\Gamma_1$ leaves the neighborhood. Likewise crossing the threshold of the third quadrant ($-a_2T_{34}(\mu_2))<\mu_3<0$) maintains the tendency of $\Gamma_1$ to the periodic orbit, while $\Gamma_2$ forms a central connection from $p_2$ to $pq_1$.

\subsection{Changes in the qualitative dynamics with one hyperbolic saddle ($\mu_1<0$)} \label{sec:homhom2}

When the parameter $\mu_1<0$ the system has undergone a saddle-node bifurcation, causing the disappearance of the saddle-node $pq_1$ along with the separatrix $\Gamma_1$. As a result the system has only one equilibrium, namely the saddle point $p_2$. Now, when varying the parameters $\mu_2$ and $\mu_3$ we only need to consider the bifurcations of the separatrix $\Gamma_2$. 
\begin{figure}[!h] 
    \centering 
    \includegraphics[width=0.99\textwidth]{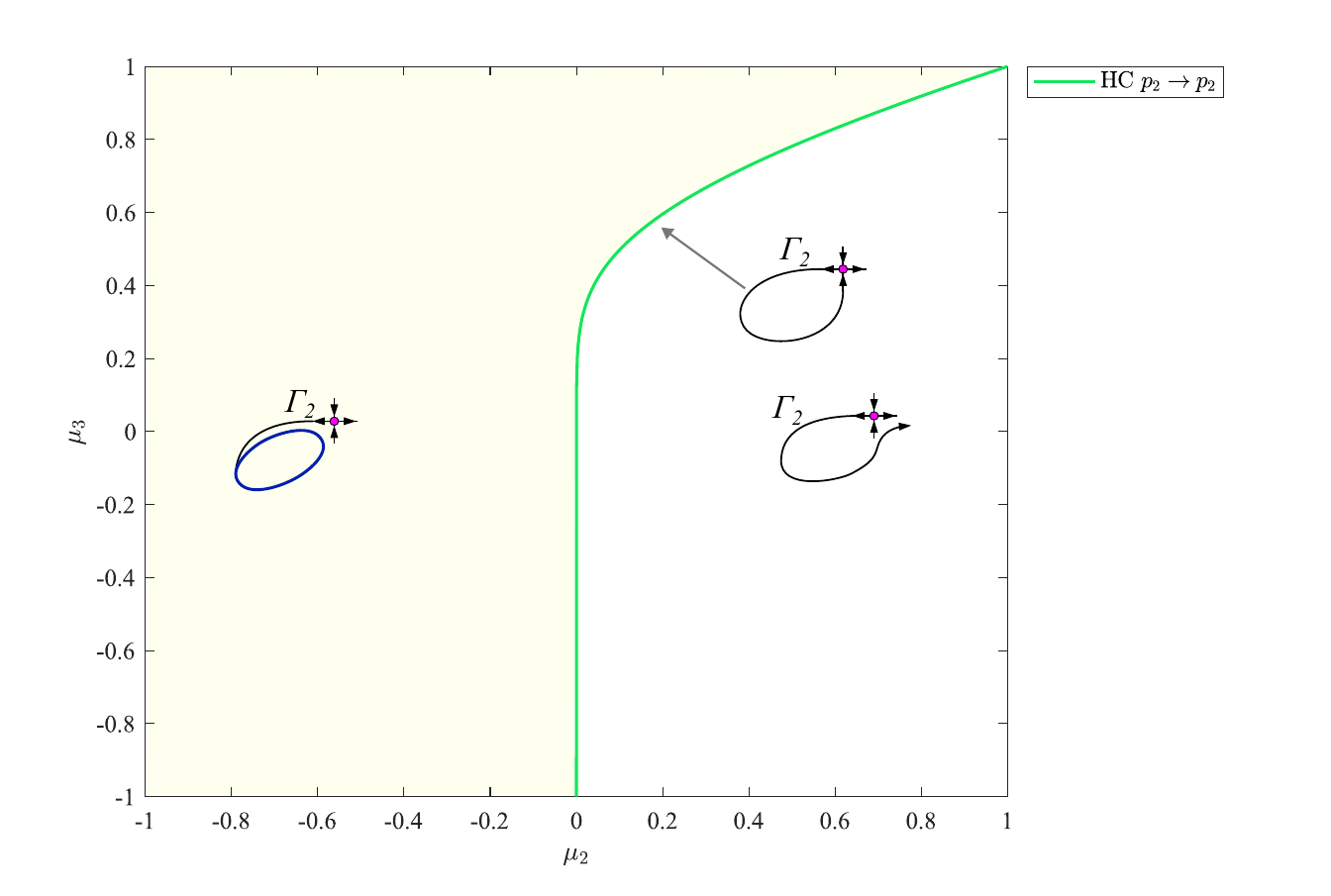} 
    \caption{Slice of the unfolding of the SNICeroclinic loop for negative $\mu_1=-0.1$ in the parameters $\mu_2$ and $\mu_3$. The green line depicts the homoclinic curve for. Left (right) to the homoclinic curve in the yellow (white) region we have persistence (destruction) of the periodic solutions. The sketches of the phase portraits illustrate the qualitative behavior of the system, depicting the behavior of the separatrices $\Gamma_1$ and $\Gamma_2$. There the pink dot identifies hyperbolic saddle equilibrium, while dark blue closed loop shows the presence of the periodic orbits.} 
    \label{fig:mu1neg}
\end{figure}

First, we consider the conditions for the persistence of the homoclinic loop with $p_2$. The conditions are similar to the ones discussed in \autoref{sec:hom1}, i.e. the fixed point of the return map $R_3$ must be zero, and hence we get the expression for the parameter $\mu_2$ in terms of $\mu_3$: 
\begin{equation}
\mu_2 = -a_1T_{12}(\mu_3). \label{cond_hc2}
\end{equation}
Unlike $\mu=0$ case, we do not need to restrict the values of $\mu_3$ as there is no other point for the separatrix $\Gamma_2$ to interact with. The connection map in \autoref{cond_hc2} depends on parameter $\mu_1$, altering the position of the homoclinic curve in the parameter space ($\mu_2$, $\mu_3$) with the changes of $\mu_1$ (\autoref{fig:mu1neg}). If $\mu_2 < -a_1T_{12}(\mu_3)$ the separatrix $\Gamma_2$ stays in the neighborhood $\mathcal{U}$ without passing through $p_2$, which identifies the presence of the periodic solutions. Otherwise, when $\mu_2 > -a_1T_{12}(\mu_3)$ the separatrix $\Gamma_2$ leaves the neighborhood $\mathcal{U}$, hence there are no invariant objects present. In summary, when $\mu_1>0$ the system can only be in three states: persistence of periodic solutions, homoclinic loop with $p_2$, which then serves as a mechanism for destruction of these solution. 

\subsection{Changes in the qualitative dynamics with two hyperbolic saddles ($\mu_1>0$)}

\begin{figure}[!h] 
    \centering 
    \includegraphics[width=0.99\textwidth]{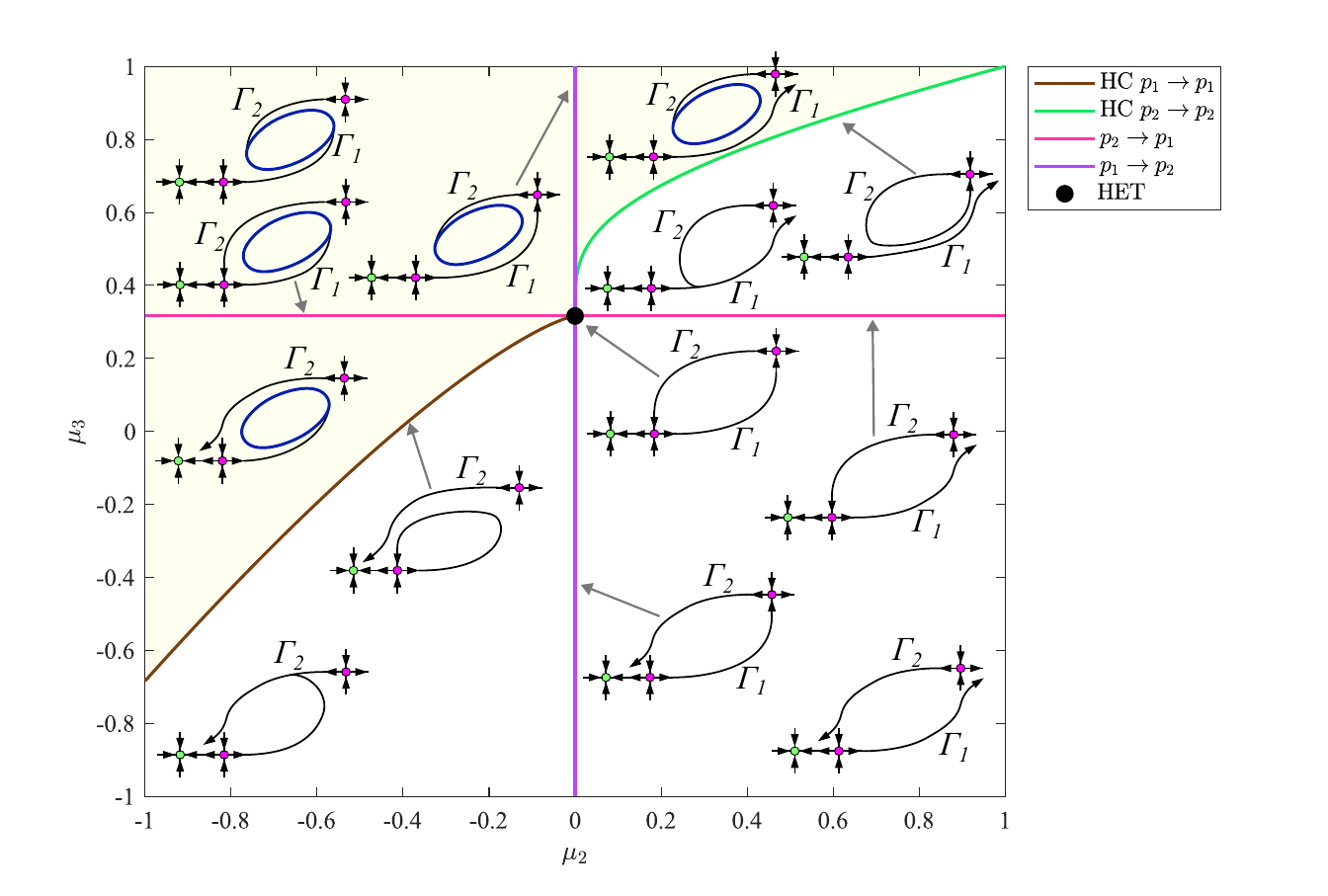} 
    \caption{Slice of the unfolding of the SNICeroclinic loop for $\mu_1>0$, when two hyperbolic (pink dots) equilibria present. The green dot identifies a stable node near saddle equilibria $p_1$. The parameter $\mu_1=0.1$ is fixed, while changes in parameters $\mu_2$ and $\mu_3$ alter the qualitative behavior of the system. Black dot signifies the heteroclinic loop between the two equilibria, the curves separate the parameter regions of different dynamics as described in the legend. Light yellow and white colored regions signify the presence and absence of a periodic orbit, respectively. The dynamics of the system is also illustrated via phase portraits, depicting the behavior of the separatrices $\Gamma_1$ and $\Gamma_2$. Dark blue closed curves identify the presence of a periodic orbit. HC: homoclinic loop, HET: heteroclinic loop.} 
    \label{fig:mu1pos}
\end{figure}

Finally, we discuss the case when $\mu_1>0$ where there are two hyperbolic equilibria $p_1$ and $p_2$. When we vary the splitting parameters $\mu_2$ and $\mu_3$ we can observe a range of different dynamic behaviors that arise from different interactions of the separatrices $\Gamma_1$ and $\Gamma_2$. We summarize all of the possible cases in the schematic in \autoref{fig:mu1pos}. In the following subsections we discuss the conditions for each of the cases.  

\subsubsection{Persistence of the heteroclinic connections}

As we fix $\mu_1>0$ the saddle-node becomes a saddle equilibrium and changes its location in the phase space. Specifically, the location of the saddle point on the cross-section $\Sigma_1$ is $x=\sqrt{\mu_1}$. To preserve the heteroclinic loop we need to move the separatrix $\Gamma_2$ so that the flow hits $\Sigma_1$ at $x=\sqrt{\mu_1}$. It can be achieved by setting the splitting parameter $\mu_3 = \sqrt{\mu_1}$. The maps $R_2$, $R_3$, $R_4$ will have a stable fixed point at $x=0$, while $R_1$'s stable fixed point will be at $x=\sqrt{\mu_1}$, indicating bidirectional heteroclinic connection. Moreover, fixing $\mu_3=\sqrt{\mu_1}$ and varying $\mu_2 \neq 0 $ we preserve the heteroclinic connection $p_2 \rightarrow p_1$, formed by the separatrix $\Gamma_2$. However, for $\mu_2<0$ the separatrix $\Gamma_1$ tends to the periodic orbit, while for for $\mu_2>0$ it leaves the neighborhood $\mathcal{U}$ (see \autoref{fig:mu1pos}). On the other hand, if we fix $\mu_2=0$ and vary the splitting parameter $\mu_3$ we only preserve the heteroclinic connection $p_1 \rightarrow p_2$, formed by the separatrix $\Gamma_1$. The separatrix $\Gamma_2$ either tends to the periodic orbit for $\mu_3>\sqrt{\mu_1}$ or leaves the neighborhood $\mathcal{U}$ for $\mu_3<\sqrt{\mu_1}$.  

\subsubsection{Homoclinic loop $p_1 \rightarrow p_1$}

As we have two saddle equilibria we have a possibility for two distinct homoclinic loop - $p_1 \rightarrow p_1$ and $p_2 \rightarrow p_2$. In the first case we need to ensure the intersection of the stable and unstable manifolds of the saddle $p_1$. For that we require that the trajectories in the neighborhood $\mathcal{U}$ are to enter the cross-section $\Sigma_1$ at a point $(\sqrt{\mu_1},\epsilon)$, where the location of the saddle equilibrium is, so the return map $R_1$ has to has a fixed point at $x=\sqrt{\mu_1}$:
\begin{align}
    0 & = R_1(\sqrt{\mu_1})-\sqrt{\mu_1} \\
    0 & = a_2 T_{34}(\mu_2) + \mu_3 -\sqrt{\mu_1} \\
    \mu_3 & = -a_2 T_{34}(\mu_2) + \sqrt{\mu_1}.
\end{align}
Additionally, we require the separatrix $\Gamma_1$ not to leave the neighborhood $\mathcal{U}$, hence $\mu_2<0$. This ensures there exists a homoclinic connection $p_1$ to $p_1$. The separatrix $\Gamma_2$ leaves the neighborhood $\mathcal{U}$. Setting $\mu_3 > -a_2 T_{34}(\mu_2) + \sqrt{\mu_1}$ leads to the existence of a periodic solution in the system. Specifically, the separatrix $\Gamma_1$ tends to a periodic orbit, while $\Gamma_2$ leaves the neighborhood as before. In the case when  $\mu_3 < -a_2 T_{34}(\mu_2) + \sqrt{\mu_1}$ the stable point in the return map $R_1$ disappears, indicating that the separatrix $\Gamma_1$ overlaps with $\Gamma_2$ and leaves the neighborhood $\mathcal{U}$. Thus the surface $\mu_3 < -a_2 T_{34}(\mu_2) + \sqrt{\mu_1}$ indicated a mechanism of birth and destruction of periodic orbits in the system. 

\subsubsection{Homoclinic loop $p_2 \rightarrow p_2$}

The conditions for the second homoclinic loop $p_2 \rightarrow p_2$ require the intersection of the stable and unstable manifolds of the saddle point $p_2$, hence the return map $R_3$ needs to have a stable point at $x=0$, as was discussed in \autoref{sec:hom1} and \autoref{sec:homhom2}: 
\begin{equation}
    \mu_2  = -a_1T_{12}(\mu_3).
\end{equation}
Moreover, we need to impose a condition for $\Gamma_2$ not to leave the neighborhood $\mathcal{U}$, thus $\mu_3>\sqrt{\mu_1}$. Under these conditions the separatrix $\Gamma_2$ forms a homoclinic loop with $p_2$, and the separatrix $\Gamma_1$ leaves the neighborhood $\mathcal{U}$. If we set $\mu_2  < -a_1T_{12}(\mu_3)$ the separatrix $\Gamma_2$ remains in the neighborhood $\mathcal{U}$ but does not go though $p_2$, hence it tends to a periodic orbit. Otherwise if $\mu_2  > -a_1T_{12}(\mu_3)$, the return map $R_3$ does not have a stable fixed point, indicating that the separatrix $\Gamma_2$ overlaps with $\Gamma_1$ and leaves the neighborhood $\mathcal{U}$. Here surface $\mu_2  = -a_1T_{12}(\mu_3)$ acts as a threshold for a transition from periodic solutions to their disappearance via a homoclinic loop with $p_2$. 

\subsubsection{Periodic solutions: their birth and destruction}

After considering the persistence of homoclinic loops with $p_1$ and $p_2$ as well as the heteroclinic loop $p_1 \leftrightarrow p_2$ in the system, we identify the following mechanisms governing the birth and destruction of the periodic solutions in the system for $\mu_1>0$. For $\mu_1>0$ the periodic orbits can be born/destroyed via one of the following mechanisms: homoclinic bifurcation with $p_1$, homoclinic bifurcation with $p_2$, heteroclinic connection $p_1 \leftrightarrow p_2$. Depending on the transition mechanism (whether it is a homoclinic bifurcation with $p_1$ or $p_2$), the flow, that originates in the neighborhood $\mathcal{U}$, could leave $\mathcal{U}$ from different sides. We remark that the above holds for Type I (and on time-reversal Type IV) but for Types II and III saddle-nodes bifurcations of periodic orbits (SNIPER) will appear \cite{dumortier_elementary_1994}, as in this case the periodic orbits emerging at the homoclinic bifurcations with $p_1$ and $p_2$ will be of different stabilities. Specifically, one will give rise to a stable limit cycle, and the other one will generate an unstable one. Moreover, if the parameters are such that $\mu_2<0$ and $\mu_3>\sqrt{\mu_1}$, both separatrices will tend to the limit cycle. On the other hand, if $\mu_2>0$ and $\mu_3<\sqrt{\mu_1}$, both separatrices will leave the neighborhood $\mathcal{U}$. Taken together, this demonstrates various global mechanisms responsible for birth and destruction of periodic solutions. Specifically, changes in the parameters $\mu_2$ and $\mu_3$ can lead to the convergence of the solutions to various attractors under perturbations to the limit cycle. 

\subsection{Unfolding the SNICeroclinic bifurcation in three parameters}
 \begin{figure}[!h] 
    \centering 
    \includegraphics[width=0.99\textwidth]{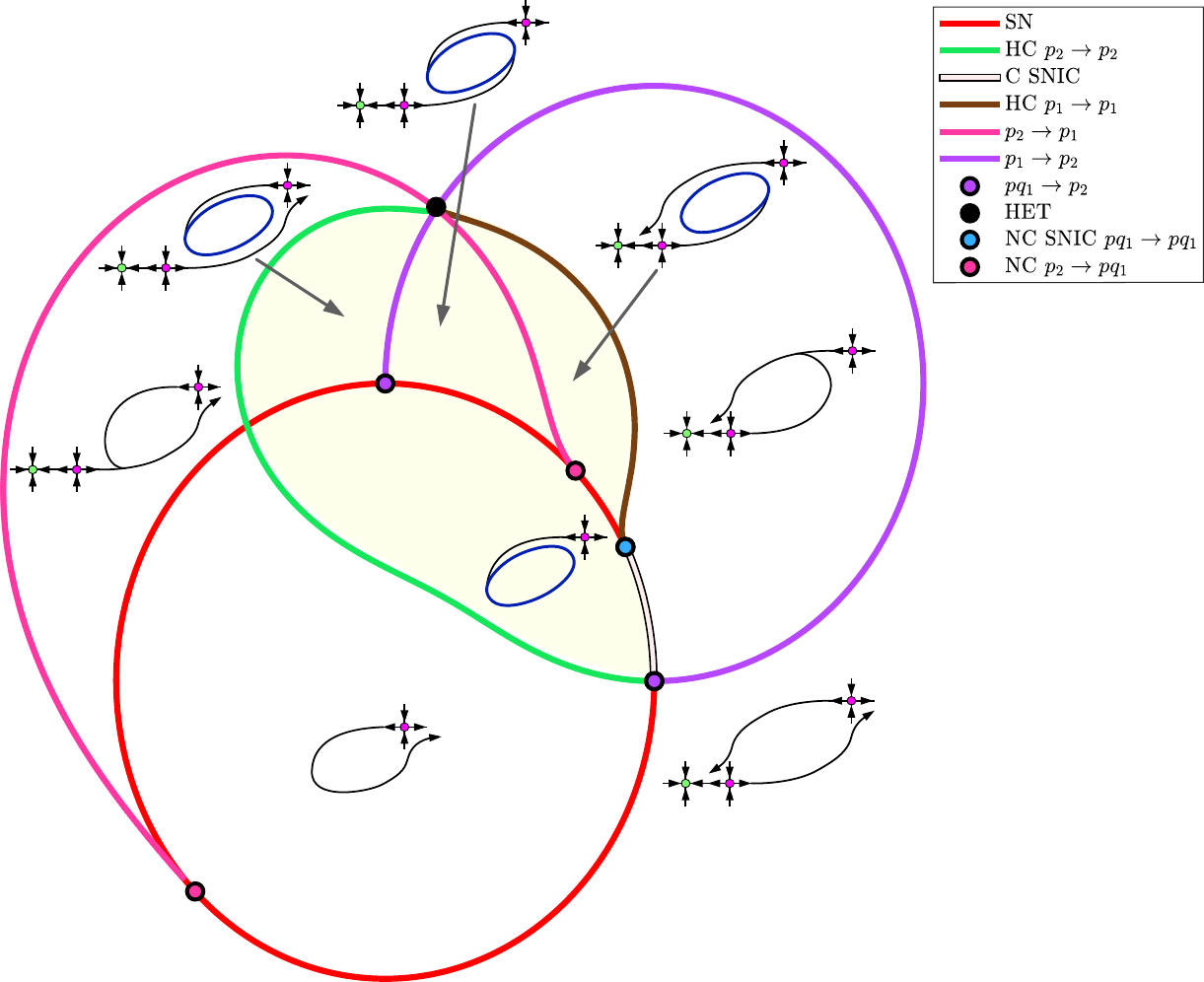} 
    \caption{Summary diagram showing intersection of $(\mu_1,\mu_2,\mu_3)$ bifurcation diagram on the boundary of a spherical neighborhood of the non-central SNICerocinic at $(0,0,0)$. This summarizes the bifurcations shown in Figure~\ref{fig:mu1neut}, \ref{fig:mu1neg} and \ref{fig:mu1pos}. One point is removed and the sphere is projected onto the plane.  The red circle signifies saddle-node bifurcation (SN) and separates region for $\mu_1<0$ (inside the red  circle, \autoref{fig:mu1neg}) and region for $\mu_1>0$ (outside the red circle, \autoref{fig:mu1pos}). Light-red region of the saddle-node circle corresponds to the central SNIC (C SNIC) with $pq_1$. Light yellow region signifies where the oscillatory behavior is observed. HC: homoclinic loop, HET: heteroclinic loop, NC: non-central.} 
    \label{fig:9}
\end{figure}
Having considered a neighborhood of the non-central SNICeroclinic bifurcation by looking at the cases for $\mu_1=0$, $\mu_1<0$, $\mu_1>0$, we superimpose our findings to obtain a codimension-three unfolding (\autoref{fig:9}). The projection highlights the transition between cases from $\mu_1<0$ (inside the red circle) to $\mu=0$ (on the red circle boundary) to $\mu_1>0$ (outside the red circle). Specifically, we observe the persistence of the region with limit cycle solutions in the 3-parameter space (light yellow shading in (\autoref{fig:9}) and how heteroclinic, homoclinic and non-central and central SNIC give rise to oscillatory dynamics. We note that as we have a transition from $\mu_1=0$ to $\mu_1>0$ the non-central SNICeroclinic loop transforms into a standard heteroclinic loop between two saddles (black dot). This mirrors, in higher codimension, the way a codimension-two SNHL point representing non-central SNIC transitions into a central SNIC bifurcation curve \cite{schecter_saddle-node_1987}. 
Additionally, in \cite{dumortier_elementary_1994}, where the mixed case of the SNICeroclinic loop is considered, the additional SNIPER bifurcation curve is identified, adding two distinct regions compared to the diagram shown in \autoref{fig:9} and reorganising the transition between various dynamic states.

\section{Discussion}
\label{sec:discussion}

In this paper we propose definitions for (non-central) SNICeroclinic loops and present a generic unfolding of the non-central SNICeroclinic loop in the plane. This is a codimension-three bifurcation of a heteroclinic loop between a saddle and a saddle-node. Our analysis augments previous work on the saddle-node homoclinic loop unfolding \cite{schecter_saddle-node_1987,chow_bifurcation_1990,deng_homoclinic_1989}, that governs the transition from a homoclinic-born and -destroyed oscillations to oscillations born and destroyed via a SNIC. Our work parallels that in \cite{grozovski_bifurcations_1996,dumortier_elementary_1994}. We highlight that the unfolding in this scenario includes saddle-node local bifurcation as well as heteroclinic, homoclinic and SNIC global bifurcations.

We suggest a classification of different types of planar non-central SNICeroclinic, while our focus of unfolding is on Type I. This is the simplest and most stable case (the saddle-node has linearly stable eigenvalue and the saddle is dominantly stable), and we give an example of Type I in \autoref{sec:examples}. We also provide an example in \autoref{sec:examples} of Type IV that we identified in a GTPase activation model \cite{zmurchok_coupling_2018,zmurchok_local_2023}. The unfolding of mixed cases Type II and III can be found in the work of \cite{dumortier_elementary_1994} who show that, in addition to the unfolding structure that we have identified in Type I, Type II and III also feature a SNIPER (saddle-node of periodic orbits) curve, where stable an unstable periodic orbits meet. The SNIPER bifurcation has been identified as a route to hidden attractors both theoretically and experimentally in electronic circuits \cite{kumarasamy_saddle-node_2023}. Its interaction with the SNICeroclinic bifurcation suggests new dynamical structures that may be interesting from both theoretical and experimental perspectives. 

\subsection{Extensions to higher dimensions}

The minimal dimension that permits a SNICeroclinic loop is two. As with Schechter's work on the SNHL in the planar case \cite{schecter_saddle-node_1987} that was extended to higher dimensions in \cite{chow_bifurcation_1990,deng_homoclinic_1989}, we highlight that there is a significant task to extend the analysis here to unfolding general non-central SNICeroclinic loops in a higher-dimensional system. 

The unfolding of the heteroclinic loop between non-hyperbolic and hyperbolic saddle equilibria, where the non-hyperbolic equilibrium is assumed to undergo a transcritical bifurcation, for a 4-dimensional system is presented in \cite{geng_bifurcations_2008}. Further work considers the case of non-generic heteroclinic loop with non-hyperbolic equilibria \cite{geng_bifurcations_2012}. Some other methods unfolding heteroclinic loops in higher-dimensional systems can be found in \cite{shashkov_complex_1996,homburg_chapter_2010,zhu_bifurcations_1998,jin_bifurcations_2003}. 

For $n>2$ there will be many types of SNICeroclinic loop in \autoref{eq:systemRn} that need to be considered for a complete classification. In particular, note 
{\bf H1*} means that $pq_1$ has either a center-stable, or a center-unstable, manifold of dimension $n$. The saddle $p_2$ also has several choices for invariant manifold dimension
$$
\dim W^{s}(p_2)=k_{s}\geq 1,~~\dim W^{u}(p_2)=k_{u}\geq 1,
$$
such that $k_s+k_u=n$. The simplest cases (parallel to Types I and IV in Table~\ref{tab:cases}) will be where $W^{cs}(pq_1)=n$, $k_u=1$, $\Gamma_2$ approaches in non-central direction (Type I/II) and where $W^{cu}(pq_1)=2$, $k_s=1$, $\Gamma_1$ leaves in non-central direction (Type III/IV), and this will also be present at codimension three. Even in these cases the unfolding may be richer due for example to strong and weak directions of attraction/repulsion, and rotations due to complex eigenvalues. More general cases (in particular hetero-dimensional cases) will typically have higher codimension but the presence of multiple connections will give much richer dynamical and bifurcation behavior, as observed for heteroclinic cycles between saddles \cite{bonatti_robust_2008}.

\subsection{Applications}

Our recent analysis of a mathematical model of neural circuitry in the amygdala revealed a heteroclinic loop between saddle-node and saddle \cite{nechyporenko_neuronal_2024}. This loop acts as a transition mechanism between a SNIC and a homoclinic orbit with distinct equilibria.

The analysis in the present paper has identified a range of qualitative dynamical behaviors, including central and non-central SNIC, periodic and homoclinic orbits as well as absence of invariant objects in the system. All of these behaviors stem from the non-central SNICeroclinic loop, highlighting its importance for induction of oscillatory dynamics.  Specifically, non-central SNIC is exactly a SNHL bifurcation, that has been dubbed as the organizing center for multi-pulse excitability in lasers \cite{krauskopf_excitability_2003, wieczorek_bifurcations_2005}.  

Just as in lasers, for neurons the route to oscillations could plays a critical role in understanding their experimentally observed behavior. Distinct global bifurcation mechanisms, such as homoclinic and SNIC, give rise to different excitability classes, each with specific dynamical and physiological signatures \cite{rinzel_analysis_1989,izhikevich_neural_2000}. While Type I (SNIC) and Type II (Hopf) excitability have been extensively discussed, recently it has been shown that changes in temperature induce switches from homoclinic (Type III) to SNIC action potential via SNHL \cite{hesse_temperature_2022}, underscoring the importance of global bifurcations for understanding experimental data. In our analysis, we reveal that such transition can occur via SNICeroclinic bifurcation, opening a new perspective on how oscillatory regimes emerge and disappear in neuronal dynamics.

In the polynomial system (\ref{eq:poly}) we have also demonstrated that timescale separation plays a role in inducing SNICeroclinic loop. Similar to how strong timescale separation has been shown to generate SNIC (Type I) excitability \cite{de_maesschalck_neural_2015}, our results demonstrate that reducing the degree of timescale separation leads to a transition from a non-central to a central connection at the saddle node. This structural change modifies the spike waveform, thereby providing a potential signature of the underlying fast–slow dynamics observable in neuronal recordings.  

Furthermore, perturbations to an excitable 1D medium have been shown to give rise to traveling pulses and the characteristics of their behavior to be mediated by the type of global bifurcation associated with oscillations \cite{moreno-spiegelberg_bifurcation_2022}. Also, it has previously been shown that a robust heteroclinic cycle can give rise to the family of periodic traveling waves, which is considered in Rock-Paper-Scissors model \cite{hasan_spatiotemporal_2021}. Extending research on dynamic transitions associated with SNICeroclinic loops can provide deeper theoretical insights into the mechanisms underlying of excitability and traveling wave phenomena.

A type of non-central heteroclinic between a saddle and a saddle-node, was previously described in the literature \cite{dumortier_elementary_1994}. In this study, we reveal that this results in a rich variety of bifurcations of oscillations and speculate that including the SNICeroclinic in continuation toolboxes will be a practically useful addition to the various hyperbolic homoclinic and heteroclinic bifurcations that can be found, for example, using \cite{veltz_bifurcationkitjl_2020, doedel_auto-07p_2007, de_witte_interactive_2012}. In particular, non-hyperbolic orbit continuation is a challenging task. Therefore, robust methods to numerically identify and continue the SNICeroclinic in continuation routines will be especially useful.

\section*{Data accessibility}
The code to reproduce the analysis can be found in a \href{https://github.com/ktrnnchprnk/SNICeroclinic.git}{github repository}


\section*{Declaration of AI use}
We have not used AI-assisted technologies in creating this article.

\section*{Conflict of interest declaration}

The authors declare that they have no conflicts of interest.

\section*{Disclaimer}
For the purpose of open access, the author has applied a ‘Creative Commons Attribution (CC BY) licence to any Author Accepted Manuscript version arising.

\section*{Acknowledgements}
We would like to acknowledge helpful comments and suggestions by Sebastian Wieczorek, Jan Sieber, and Claire Postlethwaite. 

\section*{Funding}
KTA and PA gratefully acknowledge the financial support of the EPSRC via grant EP/T017856/1. KN is a PhD student funded by the EPSRC.

\clearpage
\printbibliography
\end{document}